\documentclass[psamsfonts,reqno]{amsart}

\usepackage{amssymb}
\usepackage{amscd}
\usepackage{eucal}
\usepackage{psfig}
\usepackage[dvips]{graphicx}
\usepackage[small]{caption}

\newcommand{\beq}{\begin{equation}}
\newcommand{\eeq}{\end{equation}}

\newcommand{\Real}{\mbox{Re\,}}

\newcommand{\R}{\mathbb{R}}

\newcommand{\C}{\mathbb{C}}
\newcommand{\N}{\mathbb{N}}

\newcommand{\cF}{\mathcal{F}}

\newcommand{\cH}{\mathcal{H}}

\newcommand{\sfrac}[2]{{\vphantom1\smash{\lower.5ex\hbox{\small$#1$}}\over
        \vphantom1\smash{\raise.4ex\hbox{\small$#2$}}}} % alternate fraction

\newtheorem{theorem}{Theorem}
\newtheorem*{theorem-nonumber}{Theorem}
\newtheorem{proposition}{Proposition}
\newtheorem{lemma}{Lemma}

\theoremstyle{remark}

\theoremstyle{definition}

\newcommand{\rmd}{d}
\newcommand{\rmi}{i}

\begin{document}

\title[Watson Resummation of Power Series]{Watson resummation of a class of Hausdorff--transformed power series}

\author[E. De Micheli]{E. ~De Micheli}
\address[E. De Micheli]{IBF -- Consiglio Nazionale delle Ricerche \\ Via De Marini, 6 - 16149 Genova, Italy}
\email[E.~De Micheli]{demicheli@ge.cnr.it}

\author[G.A. Viano]{G. A. ~Viano}
\address[G.A. ~Viano]{Dipartimento di Fisica - Universit\`a di Genova,
Istituto Nazionale di Fisica Nucleare - sez. di Genova \\
Via Dodecaneso, 33 - 16146 Genova, Italy}
\email[G.A.~Viano]{viano@ge.infn.it}

\begin{abstract}
In this paper we study a class of Hausdorff--transformed power series whose
convergence is extremely slow for large values of the argument. We perform a
Watson--type resummation of these expansions, and obtain, by the use of the
Pollaczek polynomials, a new representation whose convergence is much faster.
We can thus propose a new algorithm for the numerical evaluation of these
expansions, which include series playing a relevant role in the computation
of the partition function in statistical mechanics. By the same procedure we
obtain also a solution of the classical Hausdorff moment problem.
\end{abstract}

\maketitle

\section{Introduction}
\label{se:introduction}
Following Fuchs, Rogosinski \cite{Fuchs,Rogosinski} and Hardy \cite{Hardy}
we say that $f(x)$ is summable to the value $S$ by the \emph{continuous Hausdorff
method} (or $f(x)$ is $\cH$--summable to $S$) if
\beq
\label{1}
g(x) = \int _0^1 f(xt)\,\rmd\chi(t) \qquad (x>0)
\eeq
tends to $S$ as $x\rightarrow\infty$. Integral (\ref{1}) is a Lebesgue--Stieltjes
integral, $f(x)$ is Borel summable and bounded in every finite interval $[0,x]$,
and $\chi(t)$ is of bounded variation in $[0,1]$. Next, in Hardy \cite{Hardy}
the following theorem is proved.

\begin{theorem}[Hardy]
\label{the:1}
In order that the transformation $(\ref{1})$ should be regular, i.e., that
$f(x)\rightarrow S$ $($for $x\rightarrow\infty)$ should imply $g(x)\rightarrow S$
$($for $x\rightarrow\infty)$, it is necessary and sufficient that
$\chi(0^+)=\chi(0)=0$ and $\chi(1) = 1$.
\end{theorem}

Let us now suppose that $f(x)$ is a function of $x$ regular on the positive real
axis, and so expressible in the form
\beq
\label{2}
f(x) = \sum_{n=0}^\infty \frac{f^{(n)}(0)}{n!}\,x^n = \sum_{n=0}^\infty a_n x^n.
\eeq
Then, we substitute expansion (\ref{2}) in integral (\ref{1}) and, in view of the
uniform convergence of the Taylor series, we can exchange the sum with the integral
and obtain
\beq
\label{3}
g(x)= \int _0^1 \sum_{n=0}^\infty a_n (xt)^n\,\rmd\chi(t) =
\sum_{n=0}^\infty a_n x^n \int _0^1 t^n \,\rmd\chi(t).
\eeq
The last integral at the r.h.s. of (\ref{3}) represents the \emph{Hausdorff moment}
$\mu_n$, i.e.,
\beq
\label{4}
\mu_n = \int _0^1 t^n \,\rmd\chi(t) = \int _0^1 t^n u(t) \,\rmd t.
\eeq
In formula (\ref{4}), $\chi(t)$ is supposed to be a real function of bounded
variation in $t \in [0,1]$, and the numbers $\mu_n$ are called
\emph{moment constant}, of rank $n$, of $\chi$. If we suppose, without loss of
generality, that $\chi(0)=0$, $\chi(1)=1$, and $\chi(0^+)=\chi(0)=0$,
so that $\chi(t)$ is continuous at the origin, then $\mu_n$ is called a
\emph{regular moment constant} (see Theorem \ref{the:1} above).
Moreover, the following theorem can be proved \cite{Hardy}.

\begin{theorem}[Hardy]
\label{the:2}
Sums, differences and products of moment constants are themselves moment
constants. The product of two regular moment constants is a regular moment constant.
\end{theorem}

Two relevant examples of \emph{regular Hausdorff transformations} are the
following \cite[Theorem 200]{Hardy}:
\beq
\label{5}
\mathrm{(i)} \hspace{2.5truecm} \mu_n = \ell\int_0^1 t^n (1-t)^{(\ell-1)}\,\rmd t = {n+\ell\choose \ell}^{-1}
\qquad (\ell > 0),\hspace{1truecm}\hfill
\eeq
which corresponds to the Cesaro transformation $C^{(\ell)}$;
\beq
\label{6}
\mathrm{(ii)} \hspace{2.5truecm}
\mu_n = \frac{1}{\Gamma(\ell)} \int_0^1 t^n
\left(\log\frac{1}{t}\right)^{(\ell-1)}\,\rmd t =
\frac{1}{(n+1)^\ell} \qquad (\ell > 0), \hspace{0.1truecm}
\eeq
which corresponds to the H\"{o}lder transformation $H^{(\ell)}$.

From (\ref{3}) we are naturally led to consider expansions of the following form:
\beq
\label{7}
\sum_{n=0}^\infty \frac{f^{(n)}(0)}{n!}\,\mu_n x^n,
\eeq
where the terms $\mu_n$ are Hausdorff moments.
If we suppose that $f^{(n)}(0) = (-1)^n$, we obtain expansions which read:
\beq
\label{8}
g(x) = \sum_{n=0}^\infty \frac{(-1)^n}{n!} \,\mu_n x^n.
\eeq
These series may be slowly convergent, for values of $x$ sufficiently large.
We thus face a serious problem of numerical analysis, which is quite relevant in
view of the fact that sums like expansion (\ref{8}) occur in several problems,
including some of physical interest. For instance,
\begin{enumerate}
\item[(a)] The Laplace transform of the functions of compact support gives rise
to sums of the form (\ref{8}), if we expand in series the exponential $e^{-xt}$.
This case appears in statistical mechanics, where the partition function is the
Laplace transform of the density of states. If the latter is a function of
compact support, as in the case of harmonic crystals, then we obtain a
representation of the partition function in terms of a power series of the
type (\ref{8}) \cite{Scalas}.
\item[(b)] Confluent hypergeometric function of the following type:
\beq
\label{9}
\Phi(1,\ell+1;-x) =
\sum_{n=0}^\infty \frac{(-1)^n}{n!} {n+\ell\choose \ell}^{-1} x^n \qquad (\ell > 0),
\eeq
are expansions of the form (\ref{8}) since the terms
$\mu_n={n+\ell\choose \ell}^{-1}$ form a Hausdorff sequence (see (\ref{5})).
\item[(c)] Hausdorff sequences can be constructed as follows \cite{Widder}.
Consider a sequence $\{\mu_n\}_0^\infty$ of (real) numbers, and denote by
$\Delta$ the forward difference operator: $\Delta\mu_n = \mu_{n+1} - \mu_n$.
Then we have
\beq
\label{10}
\Delta^k \mu_n=
\underbrace{\Delta\times\Delta\times\cdots\times\Delta}_{k-\mathop{\mathrm{times}}} \mu_n =
\sum_{m=0}^k (-1)^m {k\choose m} \mu_{n+k-m} \qquad (k=0,1,2,\ldots),
\eeq
$\Delta^0$ is the identity operator, by definition.
Now, suppose that there exists a positive constant $M$ such that
\beq
\label{11}
(n+1)^{(p-1)}\sum_{i=0}^n\left|{n\choose i}(-1)^{n-i}\Delta^{n-i}\mu_i\right|^p< M
\qquad (n=0,1,2,\ldots; p > 1).
\eeq
It can be proved \cite{Widder} that condition (\ref{11}) is necessary and
sufficient to represent the sequence $\left\{\mu_n\right\}_0^\infty$ as follows:
$\mu_n = \int_0^1 t^n u(t) \, \rmd t$ (see formula (\ref{4})),
where $u \in L^p [0,1]$. Thus we can say that the set $\{\mu_n\}_0^\infty$,
constrained by the condition (\ref{11}), forms a Hausdorff sequence.
\end{enumerate}

The main purpose of the present paper consists in performing a Watson--type
resummation of expansions of type (\ref{8}), where the set of numbers
$\{\mu_n\}_0^\infty$ is assumed to be a Hausdorff sequence generated
by a function $u(t)$ (see formulae (\ref{4}) and (\ref{11})) which belongs to
$L^{(2+\epsilon)}[0,1]$ $(\epsilon >0)$. In this case we can regard the sequence
$\{\mu_n\}_0^\infty$ as the restriction to the integers of a function
$\mu(z)$ $(z \in \C)$, which belongs to the Hardy space $H^2(\C_{-1/2})$, and
which is the unique Carlsonian interpolation \cite{Boas} of the numbers
$\{\mu_n\}_0^\infty$. We can thus perform the Watson--type resummation of
expansion (\ref{8}), and finally obtain another representation whose numerical
handling is much more convenient and effective.

The paper is organized as follows. In Section \ref{se:interpolation} we study
the Carlsonian interpolation of the Hausdorff moments $\{\mu_n\}_0^\infty$, and
expand the function $\mu(\rmi y-1/2)$ in terms of the so--called Pollaczek
functions. In Section \ref{se:watson} we perform a Watson--type resummation of
expansion (\ref{8}). In Section \ref{se:truncation} we study an appropriate
truncation procedure of the new representation obtained in
Section \ref{se:watson}.
In Section \ref{se:connection} we solve the Hausdorff moment problem by the use
of the Pollaczek polynomials \cite{Bateman,Szego}, and show the connection between
this problem and the Watson resummation of expansion (\ref{8}).
Finally, Section \ref{se:numerical} is devoted to numerical analysis and examples.

\section{Interpolation of Hausdorff moments and Hardy spaces}
\label{se:interpolation}
We prove the following theorem.

\begin{theorem}
\label{the:3}
Let the sequence $\left\{\mu_n\right\}_0^\infty$ satisfy condition $(\ref{11})$
with $p \geqslant 2+\epsilon$ $(\epsilon > 0)$. Then there exists a unique
Carlsonian interpolation of the numbers $\mu_n$, denoted by $\mu(z)$
$(z \in \C, \mu(n) = \mu_n)$, that satisfies the following conditions:
\begin{itemize}
\item[{\rm (i)}] $\mu(z)$ is holomorphic in the half--plane $\Real z > -1/2$,
continuous at $\Real z = -1/2$;
\item[{\rm (ii)}] $\mu(z)$ belongs to $L^2(-\infty,+\infty)$ for any fixed value
of $\Real z \equiv x \geqslant -1/2$;
\item[{\rm (iii)}] $\mu(z)$ tends uniformly to zero as $z$ tends to infinity
inside any fixed half--plane $\Real z \geqslant \delta > -1/2$;
\end{itemize}
\end{theorem}

\begin{proof}
If the sequence $\{\mu_n\}_0^\infty$ satisfies condition (\ref{11}) with
$p \geqslant 2+\epsilon$ $(\epsilon > 0)$, then
\beq
\label{12}
\mu_n = \int_0^1 t^n u(t) \,\rmd t,
\eeq
with $u \in L^{2+\epsilon}[0,1]$. Next, set $t = e^{-s}$ in formula (\ref{12}),
and obtain
\beq
\label{13}
\mu_n = \int_0^\infty e^{-ns} e^{-s} u(e^{-s}) \,\rmd s \qquad (n=0,1,2, \ldots).
\eeq
Therefore the numbers $\mu_n$ can be regarded as the restriction to the integers
of the following Laplace transform:
\beq
\label{14}
\mu(z) = \int_0^\infty e^{-(z+1/2)s} e^{-s/2} u(e^{-s}) \,\rmd s.
\eeq
Indeed, one has $\mu(n) = \mu_n$. By applying the Paley--Wiener theorem
\cite{Hoffman} to equality (\ref{14}), and recalling that the function
$e^{-s/2} u(e^{-s})$ belongs to $L^2[0,+\infty)$, we can conclude that
$\mu(z)$ belongs to the Hardy space $H^2(\C_{-1/2})$,
$\C_{-1/2} = \{z \in \C,\,\Real z > -1/2\}$ (see Ref. \cite{Hoffman}). We can
thus state that $\mu(z)$ is holomorphic in the half--plane $\Real z > -1/2$, and
tends uniformly to zero as $z$ tends to infinity inside any fixed half--plane
$\Real z \geqslant \delta > -1/2$. We can then apply the Carlson theorem
\cite{Boas}, and say that $\mu(z)$ is the unique Carlsonian interpolation of the
numbers $\mu_n$. Furthermore, in view of the fact that $e^{-s/2} u(e^{-s})$
belongs to $L^2[0,+\infty)$, then $\mu(-1/2+\rmi y)$ belongs to
$L^2(-\infty,+\infty)$, and, consequently, property (ii) holds true for any
fixed value of $\Real z \equiv x \geqslant -1/2$.
Finally, let us note that the function $e^{-s/2} u(e^{-s})$ belongs to
$L^1[0,+\infty)$; in fact,
$\int_0^\infty |e^{-s/2} u(e^{-s})|\,\rmd s = \int_0^1 \left|u(t)/\sqrt{t}\right|\,\rmd t < \infty$
since $u \in L^{2+\epsilon}[0,1]$ $(\epsilon > 0)$. Therefore, in view of
the Riemann--Lebesgue theorem applied to representation (\ref{14}), it follows
that the function $\mu(-1/2+\rmi y)$ $(y \in \R)$ is continuous, and thus
property (i) is proved.
\end{proof}

Let us now introduce the following set of functions:
\beq
\label{15}
\Psi_n(y) =
\frac{1}{\sqrt{\pi}}\,\Gamma\left(\frac{1}{2}+\rmi y\right) P_n^{(1/2)}(y),
\eeq
where $\Gamma(\cdot)$ denotes the Euler gamma function, and $P_n^{(1/2)}(\cdot)$
denote the Pollaczek polynomials $P_n^{(\lambda)}(\cdot)$, with $\lambda=1/2$
(see the appendix). These polynomials (in what follows the superscript
$\lambda = 1/2$ will be omitted) are orthonormal in $(-\infty,+\infty)$ with
weight function \cite{Bateman,Szego}
\beq
\label{16}
w(y) = \frac{1}{\pi} \left|\Gamma\left(\frac{1}{2}+\rmi y\right)\right|^2.
\eeq
Therefore the orthonormality condition reads:
\beq
\label{17}
\int_{-\infty}^{+\infty} w(y) P_n(y) P_m(y)\,\rmd y = \delta_{n,m}.
\eeq
It can be proved \cite{Itzykson} that the functions $\{\Psi_n(y)\}_0^\infty$
form a complete basis in the space $L^2(-\infty,+\infty)$.
Therefore the function $\mu(-1/2+\rmi y)$, which belongs to $L^2(-\infty,+\infty)$
(see Theorem \ref{the:3}), can be expanded in terms of this basis. We can state
the following proposition.

\begin{proposition}
\label{pro:1}
If the sequence $\{\mu_n\}_0^\infty$ satisfies condition $(\ref{11})$ with
$p \geqslant 2+\epsilon$ $(\epsilon>0)$, then
\beq
\label{18}
\mu\left(-\frac{1}{2}+\rmi y\right) = \sum_{n=0}^\infty c_n \Psi_n(y),
\eeq
which converges in the $L^2$--norm. The coefficients $c_n$ are given by
\beq
\label{19}
c_n =
\frac{1}{\sqrt{\pi}}\int_{-\infty}^{+\infty} \mu\left(-\frac{1}{2}+\rmi y\right)
\Gamma\left(\frac{1}{2}-\rmi y\right) P_n(y)\,\rmd y.
\eeq
\end{proposition}

\begin{proof}
The sequence $\{\mu_n\}_0^\infty$ is a Hausdorff sequence satisfying condition
(\ref{11}) with $p \geqslant 2+\epsilon$ $(\epsilon>0)$, then $\mu(-1/2+\rmi y)$
belongs to $L^2(-\infty,+\infty)$ (statement (ii) of Theorem \ref{the:3}), and
expansion (\ref{18}) converges in the sense of the $L^2$--norm; the coefficients
$c_n$ are then obtained by the use of the orthonormality condition (\ref{17}).
\end{proof}

The coefficients $c_n$ can be evaluated as follows.

\begin{theorem}
\label{the:4}
The following equality holds true:
\beq
\label{20}
c_n =
2\sqrt{\pi}\sum_{k=0}^\infty\frac{(-1)^k}{k!}\,\mu_k
P_n\left[-\rmi\left(k+\frac{1}{2}\right)\right],
\eeq
where $P_n(\cdot)$ are the Pollaczek polynomials.
\end{theorem}

\begin{figure}[tb]
\begin{center}
\includegraphics[width=12cm]{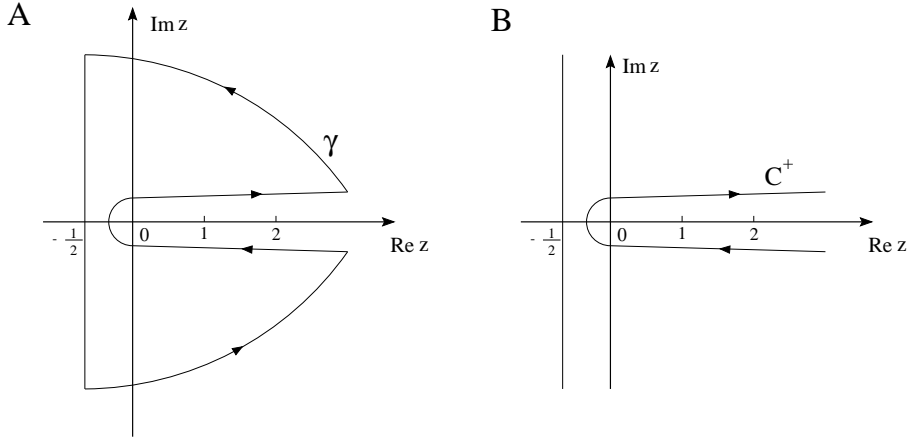}
\caption{\label{fig_1} (A) Contour integration used
for evaluating the coefficients $c_n$ (see Theorem \protect\ref{the:4}).
(B) $C^+$ is the integration path of the integral in (\protect\ref{21}).}
\end{center}
\end{figure}

\begin{proof}
Integral (\ref{19}) can be evaluated by means of the method of contour integration.
Set in formula (\ref{19}): $-1/2+\rmi y = z$, and accordingly $y = -i(z+1/2)$.
Then, performing an integration along the contour $\gamma$ shown in
Fig. \ref{fig_1}A, and taking into account the asymptotic behavior of the
gamma function, we obtain
\beq
\label{21}
\int_{-\frac{1}{2}-\rmi\infty}^{-\frac{1}{2}+\rmi\infty}
\mu(z)\Gamma(-z)P_n\left[-\rmi\left(z+\frac{1}{2}\right)\right]\,\rmd z =
\int_{C^+}\mu(z)\Gamma(-z)P_n\left[-\rmi\left(z+\frac{1}{2}\right)\right]\,\rmd z,
\eeq
where $C^+$ is a path which encircles the real positive semi--axis of the
$z$--plane in counterclockwise sense (see Fig. \ref{fig_1}B).
Then, using the theorem of residues we get
\beq
\label{22}
-\frac{\rmi}{\sqrt{\pi}}\int_{C^+}\mu(z)
\Gamma(-z)P_n\left[-\rmi\left(z+\frac{1}{2}\right)\right]\,\rmd z =
2\sqrt{\pi}\sum_{k=0}^\infty\frac{(-1)^k}{k!}\,\mu_k
P_n\left[-\rmi\left(k+\frac{1}{2}\right)\right].
\eeq
\end{proof}

\section{Watson resummation of a class of $\cH$--transformed power series}
\label{se:watson}
We prove the following theorem.

\begin{theorem}
\label{the:5}
Expansion $(\ref{8})$, where the terms $\{\mu_n\}_0^\infty$ form a Hausdorff
sequence satisfying condition $(\ref{11})$ with $p \geqslant 2+\epsilon$
$(\epsilon>0)$, can be rewritten in the following form:
\beq
\label{23}
g(x) = \sum_{n=0}^\infty\frac{(-1)^n}{n!}\mu_n x^n =
\frac{\sqrt{2}}{x+1}\sum_{n=0}^\infty u_n \rmi^n \left(\frac{x-1}{x+1}\right)^n
\qquad (x>0),
\eeq
where
\beq
\label{24}
u_n=\sqrt{2}\sum_{k=0}^\infty\frac{(-1)^k}{k!}\mu_k
P_n\left[-\rmi\left(k+\frac{1}{2}\right)\right],
\eeq
$P_n(\cdot)$ being the Pollaczek polynomials. The convergence of expansion
$(\ref{23})$ is uniform on any compact subset of the real positive axis.
\end{theorem}

\begin{proof}
We start by rewriting expansion (\ref{8}) in the following form:
\beq
\label{25}
\sum_{n=0}^\infty\frac{(-1)^n}{n!}\mu_n x^n =
\sum_{n=0}^\infty \frac{(-1)^n}{n!} \mu_n e^{n\alpha} \qquad (\alpha = \ln x).
\eeq
Next, using once again the theorem of residues, we rewrite the sum (\ref{25})
as the following integral:
\beq
\label{26}
\sum_{n=0}^\infty \frac{(-1)^n}{n!} \mu_n e^{n\alpha} =
\frac{1}{2\pi\rmi}\int_{C^+} \Gamma(-z) \mu(z) e^{\alpha z}\,\rmd z,
\eeq
where the path $C^+$ encircles the real positive semi--axis of the $z$--plane
(see Fig. \ref{fig_1}B). Equality (\ref{26}) holds true since:
\begin{itemize}
\item[{\rm (i)}] $\mu(z)$ is the Carlsonian interpolation of the
moments $\{\mu_n\}$;
\item[{\rm (ii)}] for $z = n$ ($n=0,1,2,\ldots$), one has $\mu(n)=\mu_n$,
and the function $\Gamma(-z)=\Gamma(-n)$ is singular and has simple poles with
residues $(-1)^n/n!$.
\end{itemize}
We can now close the contour $C^+$ as shown in Fig. \ref{fig_1}A. We have,
by exploiting the asymptotic behavior of the gamma function
\beq
\label{27}
\oint_\gamma \Gamma(-z) \mu(z) e^{\alpha z} \,\rmd z = 0.
\eeq
From (\ref{27}), and using the Stirling formula for the gamma function,
\beq
\label{28}
\begin{split}
\int_{C^+} \Gamma(-z) \mu(z) e^{\alpha z}\,\rmd z &=
\int_{-\frac{1}{2}-\rmi\infty}^{-\frac{1}{2}+\rmi\infty} \Gamma(-z) \mu(z)
e^{\alpha z}\,\rmd z \\
&= \rmi\int_{-\infty}^{+\infty} \Gamma\left(\frac{1}{2}-\rmi y\right)
\mu\left(\rmi y-\frac{1}{2}\right) e^{\alpha(\rmi y-1/2)}\,\rmd y.
\end{split}
\eeq
In the latter integral we use (\ref{15}) and (\ref{18}).
Then from (\ref{26}) and (\ref{28}) we have
\beq
\label{29}
\begin{split}
\sum_{n=0}^\infty \frac{(-1)^n}{n!} \mu_n e^{n\alpha} &=
\frac{1}{2\pi}\int_{-\infty}^{+\infty}\Gamma\left(\frac{1}{2}-\rmi y\right)
\mu\left(\rmi y-\frac{1}{2}\right) e^{\alpha(\rmi y-1/2)}\,\rmd y \\
&=\frac{e^{-\alpha/2}}{2\pi\sqrt{\pi}}\int_{-\infty}^{+\infty} \sum_{n=0}^\infty
c_n \Gamma\left(\frac{1}{2}-\rmi y\right) \Gamma\left(\frac{1}{2}+\rmi y\right)
P_n(y) e^{\rmi\alpha y}\,\rmd y.
\end{split}
\eeq
Using the formula
\beq
\label{30}
\Gamma\left(\frac{1}{2}-\rmi y\right) \Gamma\left(\frac{1}{2}+\rmi y\right) =
\frac{\pi}{\cosh(\pi y)},
\eeq
from (\ref{29}) it follows that
\beq
\label{31}
\int_{-\infty}^{+\infty}\sum_{n=0}^\infty c_n
\frac{P_n(y)}{\cosh(\pi y)} e^{\rmi\alpha y}\,\rmd y=
2\pi\cF^{-1}\left\{\sum_{n=0}^\infty c_n \frac{P_n(y)}{\cosh(\pi y)}\right\},
\eeq
where $\cF$ denotes the Fourier integral operator. Interchanging integration
and summation, we have, from formulae (\ref{29})--(\ref{31}):
\beq
\label{32}
\begin{split}
\frac{e^{-\alpha/2}}{2\sqrt{\pi}} \sum_{n=0}^\infty c_n
\int_{-\infty}^{+\infty}\frac{P_n(y)}{\cosh(\pi y)} e^{\rmi\alpha y}\,\rmd y &=
\frac{e^{-\alpha/2}}{2\sqrt{\pi}}\sum_{n=0}^\infty c_n
P_n\left(-\rmi\frac{\rmd}{\rmd\alpha}\right)
\left[\frac{1}{\cosh(\alpha/2)}\right] \\
&= \frac{e^{-\alpha/2}}{2\sqrt{\pi}\cosh(\alpha/2)}
\sum_{n=0}^\infty c_n \rmi^n\left[\tanh\left(\frac{\alpha}{2}\right)\right]^n.
\end{split}
\eeq
Substituting $\alpha = \ln x$ in (\ref{32}) yields
\beq
\label{33}
\sum_{n=0}^\infty \frac{(-1)^n}{n!} \mu_n x^n =
\frac{\sqrt{2}}{x+1}\sum_{n=0}^\infty u_n\rmi^n\left(\frac{x-1}{x+1}\right)^n
\qquad (x>0),
\eeq
where
\beq
\label{34}
u_n = \sqrt{2} \sum_{k=0}^\infty \frac{(-1)^k}{k!} \mu_k
P_n\left[-\rmi\left(k+\frac{1}{2}\right)\right].
\eeq
It remains to prove that the series at the r.h.s. of formula (\ref{33})
converges uniformly on any compact subset of the positive real axis.
Using the Schwarz inequality,
\beq
\label{35}
\left|\sum_{n=0}^\infty u_n\rmi^n\left(\frac{x-1}{x+1}\right)^n\right| \leqslant
\left(\sum_{n=0}^\infty |u_n|^2\right)^{1/2} \left(\sum_{n=0}^\infty
\left|\frac{x-1}{x+1}\right|^{2n}\right)^{1/2}.
\eeq
The sum $\sum_{n=0}^\infty |u_n|^2$ can be shown to be finite by using the
Parseval equality on expansion (\ref{18}). The sum
$\sum_{n=0}^\infty \left|\frac{x-1}{x+1}\right|^{2n}$
can be easily reduced to the series $\sum_{n=0}^\infty y^n$,
$y=(\frac{x-1}{x+1})^2$, which is uniformly convergent on any compact set
$y \leqslant y_0 < 1$.
\end{proof}

\section{Truncation of the resummed expansion}
\label{se:truncation}  % IV
We hereafter assume that only a finite number of Hausdorff moments $\mu_k$
(see formula (\ref{34})) are given, and, furthermore, we suppose that they can
also be affected by noise, being typically round--off numerical errors.
Accordingly, they will be denoted by $\mu_k^{(\eta)}$, $\eta$ denoting the order
of magnitude of the numerical noise. Precisely, we state:
$\left|\mu_k-\mu_k^{(\eta)}\right|\leqslant\eta$
($k=0,1,2,\ldots,k_0;\,\eta>0$); $(k_0+1)$ is the number of Hausdorff moments
which are supposed to be known. Next, we introduce the following finite sums:
\beq
\label{36}
u_n^{(\eta,k_0)}=
\sqrt{2}\sum_{k=0}^{k_0} \frac{(-1)^k}{k!} \mu_k^{(\eta)}
P_n\left[-\rmi\left(k+\frac{1}{2}\right)\right].
\eeq
With obvious notation we write: $u_n^{(0,\infty)} = u_n$.
Then, the following two auxiliary lemmas can be proved.
\begin{lemma}
\label{lem:1}
The following statements hold true: \\
%{\rm (i)}
\beq
\label{37}
\hspace{-0.0truecm}
\mathrm{(i)}\hspace{0.5truecm}
\sum_{n=0}^\infty\left|u_n^{(0,\infty)}\right|^2=
\frac{1}{\pi}\int_{-\infty}^{+\infty}
\left|\mu\left(-\frac{1}{2}+\rmi y\right)\right|^2\,\rmd y = C \qquad (C = \text{\rm const}); \hspace{0.0truecm}
\eeq
\beq
\label{38}
\hspace{-0.0truecm}\mathrm{(ii)}\hspace{4.5truecm}
\sum_{n=0}^\infty\left|u_n^{(\eta,k_0)}\right|^2 = +\infty;
\hspace{3truecm}
\eeq

\beq
\label{39}
\hspace{-0.0truecm}\mathrm{(iii)}\hspace{2.5truecm}
\lim_{\substack{k_0 \rightarrow \infty \\ \eta \rightarrow 0}}
u_n^{(\eta,k_0)} = u_n^{(0,\infty)} = u_n \qquad (n=0,1,2,\ldots);
\hspace{1truecm}
\eeq
{\rm (iv)} If $m_0(\eta,k_0)$ is defined as
\beq
\label{40}
m_0(\eta,k_0)=\max\left\{m\in\N\,:\,
\sum_{n=0}^m \left|u_n^{(\eta,k_0)}\right|^2 \leqslant C\right\},
\eeq
then
\beq
\label{41}
\lim_{\substack{k_0 \rightarrow \infty \\ \eta \rightarrow 0}}
m_0(\eta,k_0) = +\infty;
\eeq
{\rm (v)} The sum
\beq
\label{42}
M_m^{(\eta,k_0)} = \sum_{n=0}^m \left|u_n^{(\eta,k_0)}\right|^2 \qquad (m \in \N),
\eeq
satisfies the following properties:
\begin{enumerate}
\item[{\rm (a)}] It increases for increasing values of $m$;
\item[{\rm (b)}] the following relationship holds true:
\beq
\label{43}
M_m^{(\eta,k_0)} \geqslant \left|u_m^{(\eta,k_0)}\right|^2
\,\substack{\sim \\ m \rightarrow +\infty}\,
\frac{1}{(k_0!)^2}\,(2m)^{2k_0} \qquad (k_0~{\rm fixed}).
\eeq
\end{enumerate}
\end{lemma}

\begin{proof}
The proof is given, with minor modifications, in Ref. \cite{DeMicheli1}.
\end{proof}

\begin{lemma}
\label{lem:2}
The following equality holds true:
\beq
\label{44}
\lim_{\substack{k_0 \rightarrow \infty \\ \eta \rightarrow 0}}
\sum_{n=0}^{m_0(\eta,k_0)}\left| u_n^{(\eta,k_0)}-u_n^{(0,\infty)}\right|^2 = 0.
\eeq
\end{lemma}

\begin{proof}
See Ref. \cite{DeMicheli1}.
\end{proof}

Next, we introduce the following notations:
\begin{eqnarray}
\label{45}
&& g(x) \equiv g^{(0,\infty)} =
\sum_{n=0}^\infty\frac{(-1)^n}{n!}\mu_n x^n =
\frac{\sqrt{2}}{x+1} \sum_{n=0}^\infty
u_n^{(0,\infty)}\rmi^n\left(\frac{x-1}{x+1}\right)^n, \\
\label{46}
&& g^{(\eta,k_0)}=\frac{\sqrt{2}}{x+1}\sum_{n=0}^{m_0(\eta,k_0)}u_n^{(\eta,k_0)}
\rmi^n\left(\frac{x-1}{x+1}\right)^n, \\
\label{47}
&& u_n^{(\eta,k_0)} = \sum_{k=0}^{k_0} \sqrt{2}\frac{(-1)^k}{k!}\mu_k^{(\eta)}
P_n\left[-\rmi\left(k+\frac{1}{2}\right)\right].
\end{eqnarray}
Then the following is true.
\begin{theorem}
\label{the:6}
The following equality holds:
\beq
\label{48}
\lim_{\substack{k_0 \rightarrow \infty \\ \eta \rightarrow 0}}
\left| g^{(\eta,k_0)}-g^{(0,\infty)} \right| = 0.
\eeq
\end{theorem}

\begin{proof}
We have
\beq
\label{49}
\begin{split}
&\left|g^{(\eta,k_0)}-g^{(0,\infty)}\right| \\
&=\left|\frac{\sqrt{2}}{x+1}\left[
\sum_{n=0}^{m_0(\eta,k_0)}u_n^{(\eta,k_0)}\rmi^n\left(\frac{x-1}{x+1}\right)^n -
\sum_{n=0}^\infty u_n^{(0,\infty)}\rmi^n\left(\frac{x-1}{x+1}\right)^n
\right]\right| \\
&\leqslant \left|\frac{\sqrt{2}}{x+1}\right|
\left[\left|\sum_{n=m_0(\eta,k_0)+1}^\infty u_n^{(0,\infty)}
\left(\frac{x-1}{x+1}\right)^n\right|\right. \\
&\qquad +\left.\left| \sum_{n=0}^{m_0(\eta,k_0)}
\left(u_n^{(\eta,k_0)}-u_n^{(0,\infty)}\right)
\left(\frac{x-1}{x+1}\right)^n\right|\right].
\end{split}
\eeq
Now, using the Schwarz inequality,
\beq
\label{50}
\begin{split}
&\left| \sum_{n=m_0(\eta,k_0)+1}^\infty\right. u_n^{(0,\infty)} \left.
\left(\frac{x-1}{x+1}\right)^n\right| \\
&\leqslant\left(\sum_{n=m_0(\eta,k_0)+1}^\infty
\left|u_n^{(0,\infty)}\right|^2\right)^{1/2}
\left(\sum_{n=m_0(\eta,k_0)+1}^\infty\left|\frac{x-1}{x+1}\right|^{2n}\right)^{1/2}.
\end{split}
\eeq
Since
\begin{equation*}
m_0(\eta,k_0) \xrightarrow[\substack{k_0 \rightarrow \infty \\ \eta \rightarrow 0}]{} +\infty,
\qquad
\sum_{n=0}^\infty\left|u_n^{(0,\infty)}\right|^2<+\infty,
\end{equation*}
and
\begin{equation*}
\sum_{n=0}^\infty\left|(x-1)/(x+1)\right|^{2n}<+\infty \qquad \text{for} \;\; x>0
\end{equation*}
(see Lemma \ref{lem:1}), the r.h.s. of formula (\ref{50}) tends to zero as
$k_0\rightarrow+\infty$, $\eta\rightarrow 0$. Now,
\beq
\label{51}
\begin{split}
&\left|\sum_{n=0}^{m_0(\eta,k_0)} \left(u_n^{(\eta,k_0)}-u_n^{(0,\infty)}\right)
\left(\frac{x-1}{x+1}\right)^n\right| \\
& \qquad\leqslant \sum_{n=0}^{m_0(\eta,k_0)}
\left| u_n^{(\eta,k_0)}-u_n^{(0,\infty)} \right|
\left|\frac{x-1}{x+1}\right|^n \\
& \qquad\leqslant \left(\sum_{n=0}^{m_0(\eta,k_0)}
\left| u_n^{(\eta,k_0)}-u_n^{(0,\infty)}\right|^2\right)^{1/2}
\left(\sum_{n=0}^{m_0(\eta,k_0)}\left|\frac{x-1}{x+1}\right|^{2n}\right)^{1/2}.
\end{split}
\eeq
From Lemma \ref{lem:2} it follows that
\beq
\label{52}
\lim_{\substack{k_0 \rightarrow \infty \\ \eta \rightarrow 0}}
\sum_{n=0}^{m_0(\eta,k_0)}
\left| u_n^{(\eta,k_0)}-u_n^{(0,\infty)} \right|^2 = 0,
\eeq
while
$\lim_{\substack{k_0 \rightarrow \infty \\ \eta \rightarrow 0}}
\sum_{n=0}^{m_0(\eta,k_0)}\left|(x-1)/(x+1)\right|^{2n}$ is finite for $x>0$.
\end{proof}

\section{Connection with the Hausdorff moment problem}
\label{se:connection}
The classical Hausdorff moment problem can be formulated as follows
\cite{Akhiezer,Shohat}:

\vspace{0.2ex}

{\bf Problem.} Given a sequence of real numbers $\{\mu_n\}_0^\infty$,
find a function $u(t)$ such that
\beq
\label{53}
\mu_n = \int_0^1 t^n u(t)\,\rmd t~~~~~(n=0,1,2,\ldots).
\eeq

\vspace{0.2ex}

This problem is ill-posed in the sense of Hadamard \cite{Hadamard}: Suppose,
for instance, that we are looking for a solution in the space $X = L^2[0,1]$,
and assume that a solution in this space exists and is unique, but
it does not depend continuously on the data. Further, in practical cases only a
finite number of moments $\{\mu_n\}_0^N$ are known. We must then look for a
solution in a finite--dimensional subspace $X_{N+1}$ of $X$. Therefore, any
function which is orthogonal to $X_{N+1}$ cannot be recovered: the solution is
not unique. From the numerical point of view, we are led to the inversion of
matrices which are severely ill--conditioned. We shall return on these questions
later. For now we assume that a countable set of noiseless moments
$\{\mu_n\}_0^\infty$ are given, and prove the following theorem.

\begin{theorem}
\label{the:7}
Suppose that the real sequence $\{\mu_n\}_0^\infty$ of Hausdorff moments satisfy
condition $(\ref{11})$ with $p \geqslant 2+\epsilon$ $(\epsilon>0)$. Then the
function $u(t)$ can be represented by the following expansion, which converges
in the $L^2$--norm:
\beq
\label{54}
u(t)=\sum_{n=0}^\infty u_n \Phi_n(t),\qquad
u_n=\sqrt{2}\sum_{k=0}^\infty\frac{(-1)^k}{k!}\mu_k
P_n\left[-\rmi\left(k+\frac{1}{2}\right)\right],
\eeq
where
\beq
\label{55}
\Phi_n(t) = \rmi^n\sqrt{2}\,e^{-t}\,L_n(2t),
\eeq
and $P_n(\cdot)$ and $L_n(\cdot)$ are the Pollaczek and the Laguerre polynomials,
respectively.
\end{theorem}

\begin{proof}
In formula (\ref{14}) set $z = -1/2+\rmi y$. Recalling that the support
of the function $f(s)=e^{-s/2}\,u(e^{-s})$ belongs to $\R^+$,
\beq
\label{56}
\mu\left(-\frac{1}{2}+\rmi y\right)=
\int_{-\infty}^{+\infty}e^{-\rmi ys}e^{-s/2}u(e^{-s})\,\rmd s =
\cF\left\{e^{-s/2}u(e^{-s})\right\},
\eeq
where $\cF$ denotes the Fourier transform operator. Let us now return to the
expansion (\ref{18}) and to formula (\ref{15}), which gives the expression of
the functions $\Psi_n(y)$. In particular, in the integral representation of the
Euler gamma function $\Gamma(1/2+\rmi y)$: i.e,
$\Gamma(1/2+\rmi y)=\int_0^{+\infty}e^{-t}t^{(\rmi y-1/2)}\,\rmd t$,
we set $t=e^{-s}$:
\beq
\label{57}
\Gamma\left(\frac{1}{2}+\rmi y\right)=
\int_{-\infty}^{+\infty} e^{-e^{-s}} e^{-s/2} e^{-\rmi sy}\,\rmd s =
\cF\left\{e^{-e^{-s}} e^{-s/2}\right\}.
\eeq
Since the function $e^{-e^{-s}} e^{-s/2}$ belongs to the Schwartz space $S_\infty$,
\beq
\label{58}
\begin{split}
&\cF^{-1}\left\{\frac{1}{\sqrt{\pi}}\Gamma\left(\frac{1}{2}+\rmi y\right)
P_n(y)\right\} \equiv \cF^{-1}\left\{\Psi_n(y)\right\} \\
& \qquad =P_n\left(-\rmi\frac{\rmd}{\rmd s}\right)\left[\frac{1}{\sqrt{\pi}}e^{-e^{-s}}
e^{-s/2}\right],
\end{split}
\eeq
therefore, from equality (\ref{56}), and recalling once again expansion
(\ref{18}), we obtain
\beq
\label{59}
e^{-s/2}u(e^{-s})=
\sum_{n=0}^\infty c_n P_n\left(-\rmi\frac{\rmd}{\rmd s}\right)
\left[\frac{1}{\sqrt{\pi}}e^{-e^{-s}} e^{-s/2}\right].
\eeq
Reverting to the variable $t=e^{-s}$, we have:
\beq
\label{60}
u(t)=\sum_{n=0}^\infty c_n \frac{1}{\sqrt{\pi t}}
\left\{P_n\left(\rmi t\frac{\rmd}{\rmd t}\right)\left[\sqrt{t}e^{-t}\right]\right\} =
\sum_{n=0}^\infty u_n \Phi_n(t),
\eeq
where $u_n = c_n / \sqrt{2\pi}$, and the functions $\Phi_n(t)$ are given by
\beq
\label{61}
\Phi_n(t) = \sqrt{2}\frac{1}{\sqrt{t}}
P_n \left( \rmi t \frac{\rmd}{\rmd t} \right)\left[ \sqrt{t}e^{-t} \right] =
\sqrt{2}\,\rmi^n\,e^{-t}\,L_n(2t),
\eeq
where the $L_n(\cdot)$ are the Laguerre polynomials. Note that the functions
$\Phi_n(t)$ form a complete basis in $L^2(0,\infty)$, and that the convergence of
expansion (\ref{60}) is in the sense of the $L^2$--norm \cite{Viano}.
\end{proof}

We now analyze the truncation of expansion (\ref{60}). We introduce the
approximation
\beq
\label{62}
u^{(\eta,k_0)}(t)=\sum_{n=0}^{m_0(\eta,k_0)} u_n^{(\eta,k_0)} \Phi_n(t),
\eeq
where the coefficients $u_n^{(\eta,k_0)}$ are given by formula (\ref{36}),
and $m_0(\eta,k_0)$
is defined by formula (\ref{40}).

\begin{theorem}
\label{the:8}
\beq
\label{63}
\lim_{\substack{k_0 \rightarrow \infty \\ \eta \rightarrow 0}}
\left\|u^{(\eta,k_0)}-u\right\|_{L^2[0,+\infty)} = 0.
\eeq
\end{theorem}

\begin{proof}
\beq
\label{64}
\left\|u^{(\eta,k_0)}-u\right\|^2_{L^2[0,+\infty)}=
\sum_{n=m_0(\eta,k_0)+1}^\infty\left|u_n\right|^2 +
\sum_{n=0}^{m_0(\eta,k_0)}\left|u_n-u_n^{(\eta,k_0)}\right|^2.
\eeq
The statement of the theorem follows from Lemmas \ref{lem:1} and \ref{lem:2}.
\end{proof}

As already remarked, the Hausdorff moment problem, formulated as above, is
severely ill-posed. In principle one could use regularization procedures
\cite{Engl,Nashed}, among which the Tikhonov's or Tikhonov--based methods
are the most popular \cite{Groetsch,Tikhonov}. Each of these procedures consists
in restricting the class of admissible solutions to a compact subspace of the
solution space (for instance, a subspace of $X=L^2[0,1]$), by introducing
suitable bounds on the solutions. However, some problems remain, and in
particular the determination of the so--called \emph{regularization parameter},
whose optimal choice requires a precise knowledge of the majorizations on the
solutions and on the noise affecting the data.

The method presented above does not make use of any a--priori knowledge on the
solution and on the data. In several cases the truncation given by formula
(\ref{62}) can be easily determined by the properties of the
\emph{truncation number} $m_0(\eta,k_0)$, illustrated by Lemma \ref{lem:1},
and, in particular, by the statements (iv) and (v) (for numerical examples,
see the next section). Then, the statement of Theorem \ref{the:8}
guarantees the convergence of approximation (\ref{62}) to the solution, as the
number of data increases to infinity, and the noise tends to zero (see (\ref{63})).
If the number of data is too small, or the noise is too large, or, finally,
if the function $u$ to be determined is irregular (i.e., presents discontinuities
of various types), the sum $M_m^{(\eta,k_0)}$ (see (\ref{43}) in Lemma
\ref{lem:1}) can present no \emph{plateau} (see next section), and the method
cannot be used. However, this negative result still provides information:
the continuity which could at best be restored with classical regularization
procedures remains extremely weak.

Returning to the Hausdorff--transformed power series of type (\ref{8}),
it is worth remarking on the following fact: the method of resummation which we
have presented is affected by the same type of ill--posedness illustrated
above in connection with the solution of the Hausdorff moment problem.
More precisely, we face the ill--posedness connected to the reconstruction of
the function $\mu(\rmi y -1/2)$ from the sequence of Hausdorff moments
$\{\mu_k\}_0^{k_0}$. We can thus advance the following critical remark:
the method of resummation presented \emph{transforms} the type of pathology
affecting expansion (\ref{8}), i.e., slow convergence, into another type of
pathology, i.e., ill-posedness.
However, the ill-posedness of the problem, \emph{cured} by the truncation procedure
presented above, has, at least in the case of regular Hausdorff transformation,
much milder effects than the slow convergence pathology on the actual goal to be
achieved: the numerical evaluation of functions of type (\ref{8}).
In other words, \emph{the regularization of the ill-posedness cures the drawbacks
of the slow--convergence}.

\section{Numerical analysis}
\label{se:numerical}
Return to formulae (\ref{23}), (\ref{24}), and set
$x=1/r$, ($r\in[1,+\infty)$).
Then
\beq
\label{65}
~~~~~~~~~~~g\left(\frac{1}{r}\right)=\sum_{n=0}^\infty
\frac{(-1)^n}{n!}\mu_n\left(\frac{1}{c}\right)^n
\left(\frac{c}{r}\right)^n=
\frac{\sqrt{2}}{c/r+1}\sum_{n=0}^\infty u_n^{(c)} \, \rmi^n
\left(\frac{c/r-1}{c/r+1}\right)^n \qquad (c\geqslant 1),
\eeq
where
\beq
\label{66}
u_n^{(c)}=\sqrt{2}\sum_{k=0}^\infty \frac{(-1)^k}{k!}\mu_k
\left(\frac{1}{c}\right)^k P_n\left[-\rmi\left(k+\frac{1}{2}\right)\right].
\eeq
Next, we take for $c$ the value $r\in[1,+\infty)$; then the rightmost expansion
in (\ref{65}) and formula (\ref{66}) reproduce once again the original series
$\sum_{n=0}^\infty\frac{(-1)^n}{n!}\mu_n\left(1/r\right)^n$,
since $P_0(\cdot)=1$. The same type of result can be obtained when the number
of moments $\mu_k$ is finite. We can thus conclude that
the expansion that we have proposed can always be reduced to the standard one
in the interval $(0,1]$, where the original series converges rapidly.
We shall now show that the resummed expansion converges much more rapidly than
the original one for $x \gg 1$. First we discuss in detail how to manage this
new expansion numerically, and, in particular, how the truncation number
$m_0(\eta,k_0)$ can be determined.

\begin{figure}[ht]
\begin{center}
\includegraphics[width=11cm]{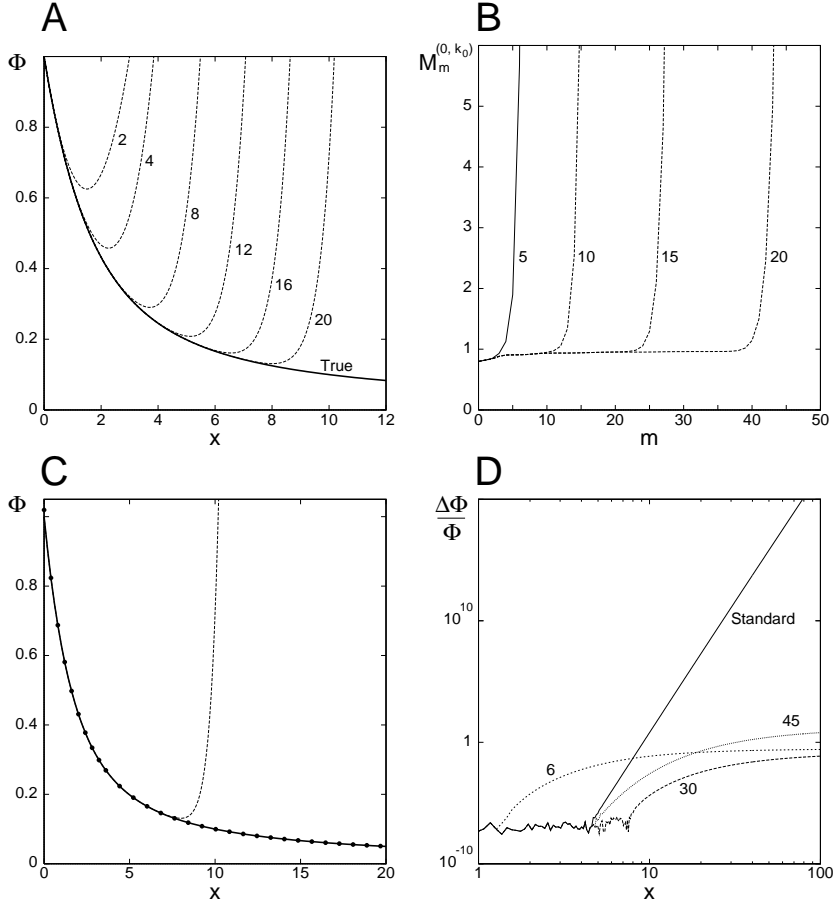}
\caption{\label{fig_2} Computation of the confluent
hypergeometric function $\Phi(1,\ell+1;-x)$ with $\ell=1$. (A)
Computations of $\Phi(1,2;-x)$ by means of expansion (\protect\ref{9})
truncated at $n=n_0$, for different values of $n_0$. The number
near each plot indicates the value of $n_0$ used for the
computation. The solid line represents the actual function
$\Phi(1,2;-x)$. (B) $M_m^{(0,k_0)}$ versus $m$ for $k_0 =
5,10,15,20$ (see (\protect\ref{42})). $(k_0+1)$ is the number of moments
$\mu_k$ used for the computation. (C) $\Phi(1,2;-x)$ computed by
means of the resummed expansion (\protect\ref{46}) (dots) by using 21
moments $\mu_k$ (i.e., $k_0 = 20$). $m_0(0,k_0)=30$ has been set
according to the analysis of the function $M_m^{(0,20)}$ in (B).
The solid line represents the actual function $\Phi(1,2;-x)$. The
dashed line shows the computation made by using expansion
(\protect\ref{9}) truncated at $n_0=20$. (D) Plot of the relative error
of reconstruction $|(\Phi_{\rm appr}-\Phi_{\rm true})/\Phi_{\rm
true}|$ in the range $x \in [1:100]$; $n_0 = k_0 = 20$. The solid
line shows the relative error computed by using expansion
(\protect\ref{9}). The dashed lines indicate the relative error of the
computations made with the resummed expansion (\protect\ref{46}), with
$m_0(0,20)=6,30,45$.}
\end{center}
\end{figure}

\begin{figure}[ht]
\begin{center}
\includegraphics[width=11cm]{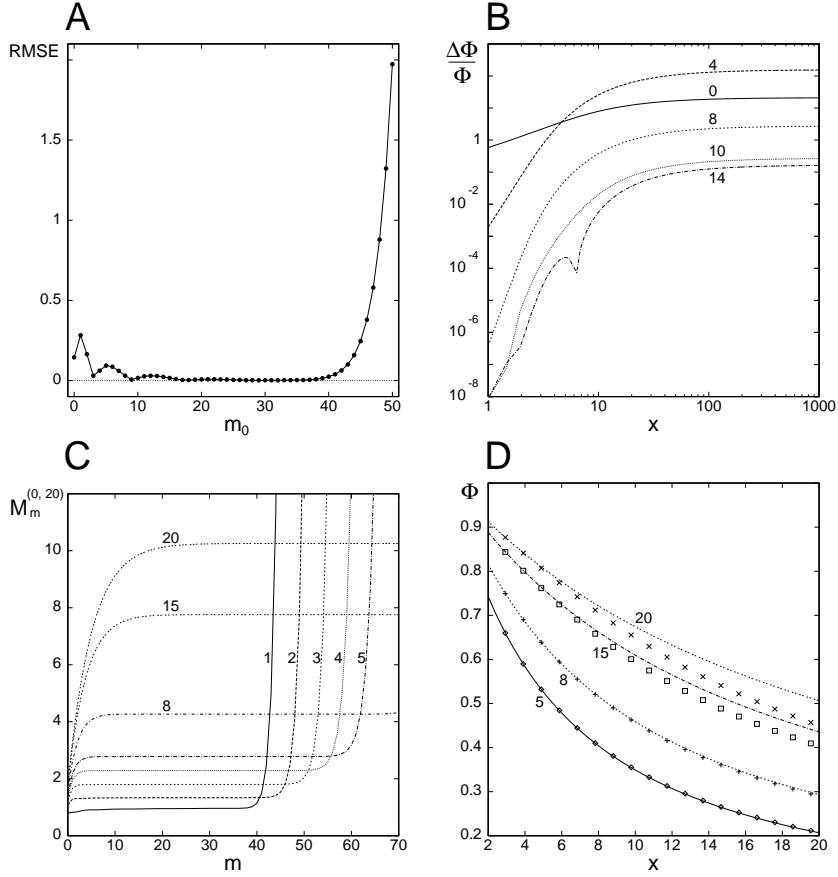}
\caption{\label{fig_3} (A) Root mean square error
of the computation made with the resummed expansion (\protect\ref{46})
with respect to the actual function $\Phi(1,2;-x)$ in the range $x
\in [0:20]$, versus $m_0(0,k_0)$; $k_0 = 20$. (B) Relative error
of computation of $\Phi(1,2;-x)$ by using the resummed expansion
(\protect\ref{46}) in the range $x \in [1:1000]$. $m_0=10$ is kept fixed,
while $k_0$ varies: $k_0=0,4,8,10,14$. (C) Computation of the
confluent hypergeometric function $\Phi(1,\ell+1;-x)$ for
different values of $\ell$. In this panel the plots of
$M_m^{(0,20)}$ (i.e., $k_0 = 20$) versus $m$ are shown. Each plot
has been computed by using the moments $\mu_k={k+\ell\choose
\ell}^{-1}$ with $\ell=1,2,3,4,5,8,15,20$. (D) Plots of
$\Phi(1,\ell+1;-x)$ computed by using the resummed expansion
(\protect\ref{46}) (dots), for various values of $\ell$: $\ell =
5,8,15,20$; $k_0 = 20$ and $m_0(0,20) = 10$ have been kept fixed.
The lines represent the actual $\Phi(1,\ell+1;-x)$.}
\end{center}
\end{figure}

\begin{figure}[ht]
\begin{center}
\includegraphics[width=11cm]{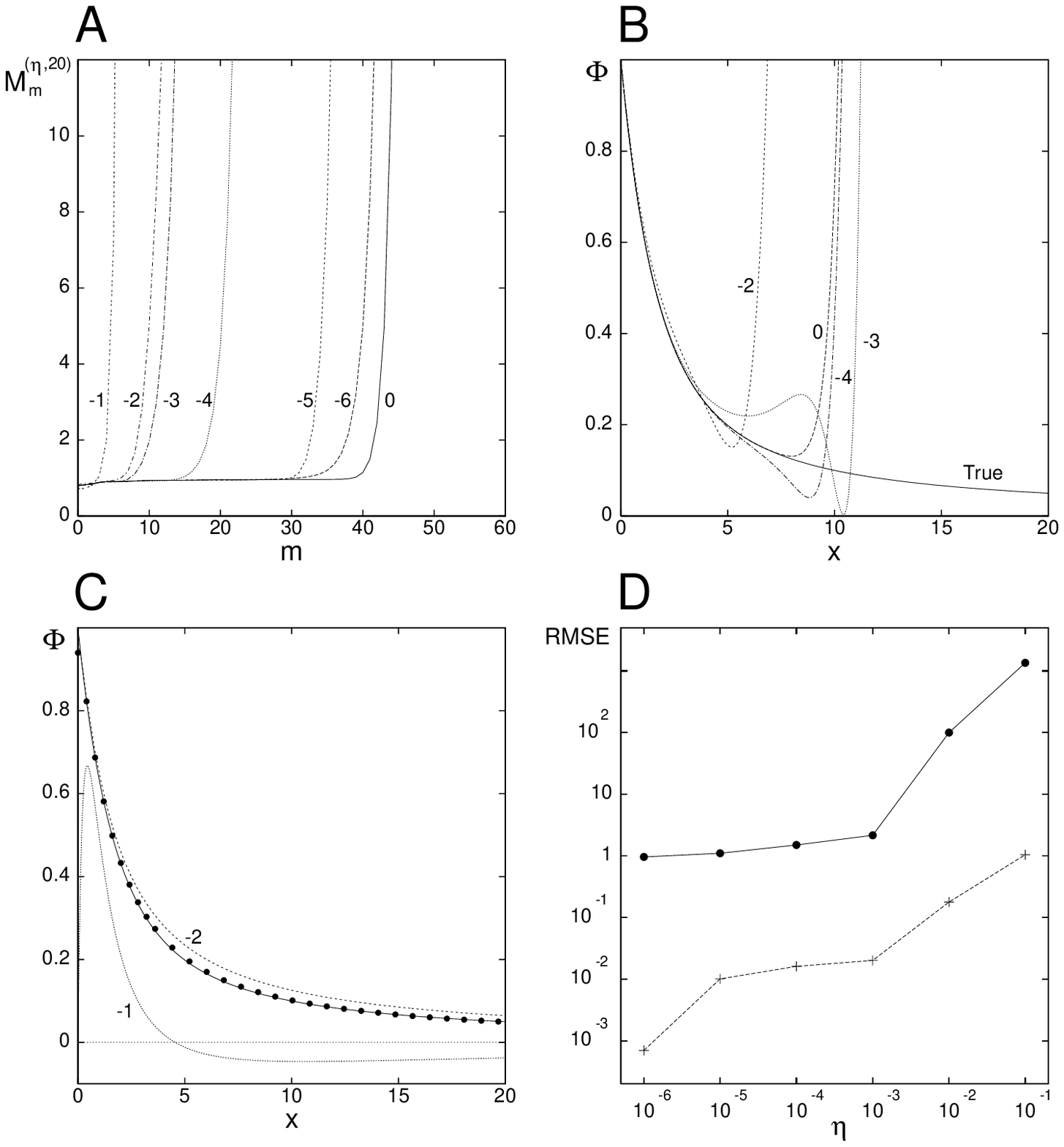}
\caption{\label{fig_4} Computation of $\Phi(1,2;-x)$
by using noisy moments $\mu^{(\eta)}_k$. The moments $\mu_k$ have
been noised by adding white noise uniformly distributed in the
interval $[-\eta,\eta]$. (A) $M_m^{(\eta,20)}$ versus $m$
computed for $\eta$ ranging from $\eta = 10^{-1}$ through $\eta =
10^{-6}$ with step $10^{-1}$. The rightmost solid line indicates
the noiseless $M_m^{(0,20)}$. (B) Comparison between the actual
function $\Phi(1,2;-x)$ (solid line) and the computations made by
using expansion (\protect\ref{9}), with $\eta = 10^{-2}, 10^{-3},
10^{-4}$. The dashed line labelled by "0" has been computed by
using the noiseless moments $\mu_k$. The number of moments used is
$(n_0+1)=21$. (C) Comparison between the actual function
$\Phi(1,2;-x)$ (solid line) and the computations made by using the
resummed expansion (\protect\ref{46}); $k_0 = 20$. The filled dots
indicate the computation made with $\eta = 10^{-3}$ and
$m_0(10^{-3},20) = 5$ (see also panel (A)). The dashed lines
represent the computations made with $\eta = 10^{-2}$,
$m_0(10^{-2},20) = 4$, and $\eta = 10^{-1}$, $m_0(10^{-1},20) =
2$. (D) Root mean square error of the computation of
$\Phi(1,2;-x)$ versus $\eta$ in the range $x \in [0:10]$; $k_0 =
n_0 = 20$. The dots indicate the error made by using expansion
(\protect\ref{9}), while the crosses represent the error made by using
expansion (\protect\ref{46}). The truncation numbers are: $m_0(10^{-6},20)
= 30$, $m_0(10^{-5},20) = 15$, $m_0(10^{-4},20) = 10$,
$m_0(10^{-3},20) = 5$, $m_0(10^{-2},20) = 4$, $m_0(10^{-1},20) =
2$.}
\end{center}
\end{figure}

The problem of evaluating $m_0(\eta,k_0)$ is intimately related to the price
that must be paid for coping with the ill--posedness of the analytic continuation
involved in the reconstruction of the function $\mu(\rmi y-1/2)$ from the
sequence of Hausdorff moments $\{\mu_k\}_0^{k_0}$.
As remarked in Section \ref{se:connection}, most regularization
procedures generally require a priori bounds on the solution and on the data.
Then one is led to introduce in the regularizing algorithm a truncation,
or a filtering, which depends on the a priori information on the solution
that one is supposed to have. Our procedure does not require any prior knowledge,
and the truncation number $m_0(\eta,k_0)$ (see Section \ref{se:truncation})
can be determined by analyzing the sum $M_m^{(\eta,k_0)}$ versus $m$
(see (\ref{42})). From statement (v) of Lemma \ref{lem:1} and from formula
(\ref{41}) it follows that, if $k_0$ is sufficiently large and $\eta$ is
sufficiently small, then $M_m^{(\eta,k_0)}$ presents a \emph{plateau}, and, after
that, it starts growing as a power of $(2k_0)$ (see (\ref{43}) and
Figs. \ref{fig_2}B, \ref{fig_3}C and \ref{fig_4}A).
The extension of the \emph{plateau} varies, and increases as the number of
moments $\mu_k$ becomes larger (see Fig. \ref{fig_2}B). Analogous results are
obtained even in the case $\eta \neq 0$ (see Fig. \ref{fig_4}A).
A simple algorithm for the automatic derivation of $m_0$ can be implemented along
the following lines. We start by observing that the knowledge of the asymptotic
behavior of $M_m^{(\eta,k_0)}$, for large $m$, allows us to restrict the range
of $m_0$ by defining an upper limit $m_\alpha$ ($m_0<m_\alpha$), which represents
the value of $m$ where approximately the asymptotic behavior sets in; in practice,
$m_\alpha$ is set as the value of $m$ where $M_m^{(\eta,k_0)}$ starts being close
enough to its asymptotic behavior. The candidate plateaux are then located by
selecting the extended intervals of $m<m_\alpha$ where the modulus of the first
numerical derivative of $M_m^{(\eta,k_0)}$ is sufficiently small.
Finally, $m_0$ is chosen as the largest value of $m$ belonging to the interval
which is closest, but inferior, to $m_\alpha$. It should be noticed that the
choice of $m_0$ within the plateau is not critical for the accuracy of the final
result. For the sake of completeness, it should also be mentioned that the
erratic behavior of the noise can produce very short plateaux located
between the true value of $m_0$ and before $M_m^{(\eta,k_0)}$ starts following
its asymptotic behavior (i.e., for $m=m_\alpha$). In this case our procedure
could fail to recover the correct value of $m_0$; this drawback has been
solved heuristically by simply rejecting plateaux shorter than a given threshold
length.

We test our method by comparing different evaluations of the confluent
hypergeometric function $\Phi(1,\ell+1;-x)$ (see formula (\ref{9})).
The latter can be directly evaluated as the Laplace transform of the
following function:
\beq
\label{67}
u(t) =
\begin{cases}
\ell(1-t)^{(\ell-1)}, & t \in [0,1], \\
0, & t \in (1,+\infty).
\end{cases}
\eeq
On the other hand, one can evaluate $\Phi(1,\ell+1;-x)$ by means of the standard
expansion at the r.h.s. of formula (\ref{9}) (truncated at a certain $n = n_0$),
and finally these results can be compared with those obtained by the Watson
resummation method (in particular, see (\ref{46})).
The results are illustrated in Figs. \ref{fig_2}--\ref{fig_4}.

In Fig. \ref{fig_2}A the plots of $\Phi(1,2;-x)$, computed by using the standard
formula (\ref{9}) with different values of $n_0$, are shown. It is evident how
the deviation from the true function rapidly explodes as $x$ increases.
Moreover, even using more moments $\mu_n$, i.e., increasing $n_0$, the
situation does not get better significantly. Figure \ref{fig_2}B shows the sum
$M_m^{(0,k_0)}$ versus $m$ for various values of $k_0$ ($k_0=5,10,15,20$).
It can be seen the presence of the \emph{plateaux}, whose length increases
as the number of moments $\mu_k$ used in the computation, i.e., $k_0$, increases
(see (\ref{41}) in Lemma \ref{lem:1}).
This figure shows how the truncation number $m_0(0,k_0)$ can be determined.
The comparison among the true $\Phi(1,2;-x)$, computed analytically as the Laplace
transform of the function $u(t)$ in (\ref{67}) (solid line), the evaluation
obtained by the truncated standard expansion (\ref{9}) (dashed line), and
truncated resummed expansion (\ref{46}) (filled dots), is shown. In this case,
$k_0=n_0=20$. From the inspection of $M_m^{(0,20)}$ in Fig. \ref{fig_2}B, it is
recovered that the plateau approximately ranges from $m=8$ through $m=38$.
The computation shown in Fig. \ref{fig_2}C was made with $m_0(0,20)=30$.
It is evident how the accuracy of the computation increases considerably when
the resummed expansion is used. This fact is made even more clear in
Fig. \ref{fig_2}D, where the relative error of computation
$\left|\Delta\Phi/\Phi\right| =
\left|(\Phi_{\rm appr}-\Phi_{\rm true})/\Phi_{\rm true}\right|$
over the range $x\in[1:100]$ is shown. The error made by using the standard
expansion (solid line) diverges, whereas the error made by using the resummed
expansion (dashed line) remains quite limited. This latter is displayed even
for the cases $m_0=6$ and $m_0=45$, which represent two values just outside the
plateau. In these cases the error increases, though in different ways, with
respect to the \emph{correct} value $m_0=30$. For any other value of $m_0$ within
the plateau, that is $m\in[8,38]$, the error of computation does not change
significantly from that at $m_0=30$. This fact is evidentiated in
Fig. \ref{fig_3}A, where the root mean square error of the computation of
$\Phi(1,2;-x)$, integrated over the interval $x\in[0,20]$, is shown for $k_0=20$,
and $m_0$ varying in the range $[0:50]$. As long as $m_0$ remains into the
plateau the error is small, whereas when $m_0$ increases far beyond the upper
limit of the plateau, the error rapidly grows.
In order to show the interplay between $k_0$ and $m_0$, in Fig. \ref{fig_3}B
it is shown the relative error for different values of $k_0$ ($k_0=0,4,8,10,14$),
and $m_0=10$ kept fixed. Even in this case the error is high when $m_0$ does not
lie within the plateau of $M_m^{(0,k_0)}$ (see also Fig. \ref{fig_2}B).
In Fig. \ref{fig_3}C the sums $M_m^{(0,20)}$ have been plotted for various
values of $\ell$ ($\ell = 1,2,3,4,5,8,15,20$), while in Fig. \ref{fig_3}D the
corresponding computations of $\Phi(1,\ell+1;-x)$ ($\ell = 5,8,15,20$), made by
means of the resummed expansion (\ref{46}) with $k_0=20$ and $m_0(0,20)=10$
are compared with the true functions. From the analysis of $M_m^{(0,20)}$ in
Fig. \ref{fig_3}C it can be seen that the value $m_0(0,20)=10$ lies outside the
plateau resulting for $\ell=15$ and $\ell=20$; correspondingly,
in Fig. \ref{fig_3}D the computation of $\Phi(1,\ell+1;-x)$ with $\ell=15$ and
$\ell=20$ clearly deviate from the corresponding true functions.
Finally, Fig. \ref{fig_4} illustrates the analysis in the case of noisy moments
$\mu_k^{(\eta)}$. To obtain the $\mu_k^{(\eta)}$, the moments $\mu_k$ have been
noised by adding white noise uniformly distributed in the interval $[-\eta,\eta]$.
Figure \ref{fig_4}A shows the sum $M_m^{(\eta,20)}$ for various values of the
noise parameter $\eta$
($\eta=10^{-1},10^{-2},10^{-3},10^{-4},10^{-5},10^{-6}$; $\ell=1$). These plots
show how the plateau gets shorter as the noise level increases.
In Fig. \ref{fig_4}B some examples of computation of $\Phi(1,2;-x)$ by using
the standard expansion for various values of $\eta$ are shown, while in
Fig. \ref{fig_4}C the computations made by using the resummed expansion are given.
This panel shows that a significant deviation from the actual function arises
only in the case of quite noisy moments ($\eta \sim 10^{-1}$),
and remains quite acceptable for small levels of noise.
The difference of accuracy achieved by using the two types of expansion are made
evident in Fig. \ref{fig_4}D, where the root mean square error, in the range
$x\in[0,10]$, is shown as a function of $\eta$.
The dots indicate the error referred to the standard expansion, while the
crosses indicate that referred to the resummed expansion. Over a wide range
of noise level, the error made by using the resummed expansion remains orders
of magnitude smaller than that made by using the standard expansion.

\renewcommand{\theequation}{A.\arabic{equation}}
\setcounter{equation}{0}
\section*{Appendix}
\label{se:appendix}
The main properties of the Pollaczek polynomials $P^{(1/2)}_n(y)$ are briefly
summarized: \\
(1) In terms of the hypergeometric series \cite{Bateman},
\beq
P^{(1/2)}_n(y)=\rmi^n ~_2F_1\left(-n,\frac{1}{2}+\rmi y,1;2\right).
\label{a1}
\eeq

\noindent
(2) They satisfy the following recurrence relation \cite{Bateman}:
\begin{subequations}
\label{a2}
\begin{eqnarray}
&&(n+1)P^{(1/2)}_{n+1}(y)-2yP^{(1/2)}_n(y)+nP^{(1/2)}_{n-1}(y)=0, \label{a2a} \\
&&P^{(1/2)}_{-1}(y)=0, \qquad P^{(1/2)}_0(y)=1.  \label{a2b}
\end{eqnarray}
\end{subequations}

\noindent
(3) The generating function is given by:
\beq
\sum_{n=0}^\infty z^n P^{(1/2)}_n(y) =
(1-\rmi z)^{(\rmi y-1/2)}(1+\rmi z)^{(-\rmi y-1/2)} \qquad (|z|<1).
\label{a3}
\eeq


\newpage

\begin{thebibliography}{20}

\bibitem{Akhiezer}
N. I. Akhiezer,
The Classical Moment Problem and Some Related Questions in Analysis,
Oliver and Boyd, Edimburgh, 1965.

\bibitem{Bateman}
A. Erd\'elyi,
Higher Trascendental Functions,
in: Bateman Manuscript Project, vol. 2,
McGraw--Hill, New York, 1954.

\bibitem{Boas}
R. P. Boas,
Entire Functions,
Academic Press, New York, 1954.

\bibitem{DeMicheli1}
E. De Micheli, G. A. Viano,
On the solution of a class of Cauchy integral equations,
J. Math. Anal. Appl. 246 (2000) 520--543.

\bibitem{Engl}
H. W. Engl,
Regularization methods for the stable solution of inverse problems,
Surveys Math. Indust. 3 (1993) 71-143.

\bibitem{Fuchs}
W. H. J. Fuchs, W. W. Rogosinski,
On typical means,
Quart. J. Math. 14 (1943) 27--48.

\bibitem{Groetsch}
C. W. Groetsch,
The Theory of Tikhonov Regularization for Fredholm Equations of the First Kind,
Pitman, Boston, 1984.

\bibitem{Hadamard}
J. Hadamard,
Lectures on the Cauchy Problem in Linear Differential Equations,
Yale Univ. Press, New Haven, CT, 1923.

\bibitem{Hardy}
G. H. Hardy,
Divergent Series,
Clarendon, Oxford, 1949,
Chapter 11 and the references quoted therein.

\bibitem{Hoffman}
K. Hoffman,
Banach Spaces of Analytic Functions,
Prentice Hall International, Englewood Cliffs, NJ, 1962.

\bibitem{Itzykson}
C. Itzykson,
Group representation in a continuous basis: An example,
J. Math. Phys. 10 (1969) 1109--1114.

\bibitem{Nashed}
M. Z. Nashed,
Generalized Inverses and Applications,
Academic Press, New York, 1976.

\bibitem{Rogosinski}
W. W. Rogosinski,
On Hausdorff's methods of summability. II,
Proc. Cambridge. Philos. Soc. 38 (1942) 344--363.

\bibitem{Scalas}
E. Scalas, G. A. Viano,
The Hausdorff moments in statistical mechanics,
J. Math. Phys. 34 (1993) 5781--5800.

\bibitem{Shohat}
J. A. Shohat, J. D. Tamarkin,
The Problem of Moments,
in: Math. Surveys Monogr., vol. 1,
American Mathematical Society, Providence, RI, 1943.

\bibitem{Szego}
G. Szeg\"{o},
Orthogonal Polynomials,
American Mathematical Society, Providence, RI, 1959.

\bibitem{Tikhonov}
A. Tikhonov, V. Arsenine,
M\'{e}thodes de R\`{e}solution de Probl\'{e}mes Mal Pos\`{e}s,
Mir, Moscow, 1976.

\bibitem{Viano}
G. A. Viano,
Solution of the The Hausdorff moment problem by the use of Pollaczek polynomials,
J. Math. Anal. Appl. 156 (1991) 410--427.

\bibitem{Widder}
D. V. Widder,
The Laplace Transform,
Princeton Univ. Press, Princeton, NJ, 1972.

\end{thebibliography}
\end{document}